\documentstyle[12pt]{article}

\newtheorem{theorem}{Theorem}

\newtheorem{definition}{Definition}

\newtheorem{corollary}{Corollary}
\newtheorem{proposition}{Proposition}

\def\EXP{{\bf E}}
\def\PROB{{\bf P}}

\def\F{{\cal F}}

\def\proof{\medskip  \par \noindent {\bf Proof.} \ \ }

\newcommand{\be}{\begin{equation}}
\newcommand{\ee}{\end{equation}}

\begin{document}
\bibliographystyle{plain}

\begin{titlepage}
\thispagestyle{empty}
\setcounter{page}{0}

\title{Queueing for ergodic arrivals and services}
\author{
L\'aszl\'o Gy\"orfi
\thanks{
Dept. of Computer Science and Information Theory,
Technical University of Budapest, 1521 Stoczek u. 2, Budapest, Hungary.
E-mail:gyorfi@szit.bme.hu}
\and
Guszt\'av Morvai
\thanks{
Research Group for Informatics and Electronics of the Hungarian Academy
of Sciences,
1521 Stoczek u. 2, Budapest, Hungary. E-mail: morvai@szit.bme.hu
}
}

\maketitle

\begin{abstract}
In this paper we revisit the results of Loynes (1962) on stability of queues
for ergodic arrivals and services, and show examples when the arrivals are 
bounded and ergodic, the service rate is constant, and under stability the
limit distribution has larger than exponential tail. 
\end{abstract}

\end{titlepage}

\section{Introduction}

The analysis of a queueing model consists of two steps: stability study 
and the characterization of the limit distribution. In this paper we consider
long range dependence, i.e. assume ergodic arrivals and services. Concerning 
stability revisit the result of Loynes (1962), who extended the Markovian 
approach in an elegant way. For the weak dependent situation the limit 
distribution usually has an almost exponential tail. Here we show 
counterexample for ergodic situation, i.e. when the arrivals are 
bounded and ergodic, the service rate is constant.

\section{Stability}

\subsection{General result}

Let $X_0$ be arbitrary random variable. 
Define
\be
X_{n+1}=(X_{n}+Z_{n+1})^+
\ee
for $n\ge 0$, where $\{Z_i\}$ is a sequence of
random variables. We are ineterested in
the stability of $\{X_i\}$, i.e. we are looking for conditions on
$\{Z_i\}$ under which $X_i$ has a limit distribution. For the classical 
Markovian approach $\{Z_i\}$ are independent and identically distributed,
when stability of $\{ X_i\}$ means that there exists a unique limit
distribution of $X_n$.  
Here consider the case when $\{Z_i\}$ is only stationary and ergodic.

Following Lindvall (1992) introduce a stronger concept of stability:
\begin{definition}
We say that the sequence $\{X_i\}$ is coupled with the sequence $\{X'_i\}$
if  $\{X'_i\}$ is stationary and ergodic
and there is an almost surely finite random variable $\tau$ such that 
\[
X'_n=X_n
\]
for $n>\tau$.
\end{definition}

Put
\begin{eqnarray*} 
V_0 & = & 0,\\
V_n & = & \sum_{i=0}^{n-1}Z_{-i},\, (n\ge 1).
\end{eqnarray*} 

\begin{theorem} \label{thm1}
If $\{Z_i\}$ is stationary and ergodic, and $\EXP\{Z_i\}<0$, then 
$\{X_i\}$ is coupled with a stationary and ergodic $\{X'_i\}$ such that 
\[
X'_0=\sup_{n\ge 0}V_n.
\]
\end{theorem}

\proof \\
{\bf Step 1.}
Let $X_{-N,-N}=0$ and define $X_{-N,n}$ for $n>-N$ by the following
recursion,
\[
X_{-N,n+1}=(X_{-N,n}+Z_{n+1})^+ \mbox{\ \ for $n\ge -N$.}
\]
We show that $X_{-N,0}$ is monoton increasing in $N$, and almost surely,
\[
\lim_{N\to\infty} X_{-N,0}=X',
\]
where 
\[
X'=\sup_{n\ge 0}V_n,
\]
and $X'$ is finite a.s.\\
Notice that 
$X_{-N,n+1}=(X_{-N,n}+Z_{n+1})^+$ for $n\ge -N$.
First we proove that for $n>-N$,  
\be
\label{expansion}
X_{-N,-n}=\max\{0,Z_n,Z_n+Z_{n-1},\dots,Z_n+\dots+Z_{-N+1}\}.
\ee
For $n=-N+1$, 
\begin{eqnarray*}
X_{-N,-N+1}&=&(X_{-N,-N}+Z_{-N+1})^+\\
&=& (Z_{-N+1})^+\\
&=&
\max\{0,Z_{-N+1}\}.
\end{eqnarray*}
For $n=-N+2$, 
\begin{eqnarray*} 
X_{-N,-N+2}&=&(X_{-N,-N+1}+Z_{-N+2})^+\\
&=& \max\{0,X_{-N,-N+1}+Z_{-N+2}\}\\
&=& \max\{0,\max(0,Z_{-N+1})+Z_{-N+2}\}\\
&=& \max\{0,Z_{-N+2},Z_{-N+2}+Z_{-N+1}\}.
\end{eqnarray*}
Now we proove by induction from $n$ to $n+1$. 
\begin{eqnarray*}
X_{-N,n+1}&=&(X_{-N,n}+Z_{n+1})^+\\
&=&\max\{0,\max\{0,Z_n,Z_n+Z_{n-1},\dots,Z_n+\dots+Z_{-N+1}\}+Z_{n+1}\}\\
&=&\max\{0,Z_{n+1},Z_{n+1}+Z_{n},\dots,Z_{n+1}+\dots+Z_{-N+1}\}.
\end{eqnarray*}
We have completed the proof of~(\ref{expansion}). Thus
\[
X_{-N,0}=\max\{0,Z_0,Z_0+Z_{-1},\dots,Z_0+\dots+Z_{-N+1}\},
\]
which imlies that $X_{-N,0}$ is monoton increasing, since the maximum 
is taken over larger and  larger set. 
It remains to prove that $X_{-N,0}$ converges to a random variable $X'$
which is finite a.s..
Now by the Birkoff strong law of large numbers for ergodic sequences, a.s. 
\[
\lim_{N\to\infty} {1\over N} \sum_{i=-N+1}^{0} Z_i=\EXP Z_1<0 
\]
(cf. Theorem 3.5.7. in Stout \cite{Sto74}), hence a.s. 
\[ 
\lim_{N\to\infty} \sum_{i=-N+1}^{0} Z_i=-\infty. 
\]
We got that there is a random variable $\tau$ such that for all
$i>0$ 
$$\infty>X_{-\tau,0}=X_{-\tau-i,0},$$
and therefore
\[
X'=\sup_{n\ge 0}V_n.
\]
{\bf Step 2.} Put
\[
X'_0=X'
\]
and for $n\ge 0$, 
\[
X'_{n+1}=(X'_{n}+Z_{n+1})^+.
\]
We show that $\{X'_i\}$ is stationary and ergodic.\\
For any sequence
$z_{-\infty}^{\infty}=(\dots,z_{-1},z_{0},z_{1},\dots)$ put
\[
F(z_{-\infty}^{\infty})
=\lim_{N\to\infty}\max\{0,z_0,z_0+z_{-1},\dots,z_0+\dots+z_{-N}\}.
\]
Then by Step 1
\[
X'_0=X'=F(Z_{-\infty}^{\infty}).
\]
We prove by induction that for $n\ge 0$,
\[
X'_n=F(T^n Z_{-\infty}^{\infty}),
\]
where $T$ is the left shift.
For $n=1$, 
\begin{eqnarray*}
F(T Z_{-\infty}^{\infty}) &=& 
\lim_{N\to\infty}
\max\{0, Z_1,Z_1+Z_{0},\dots,Z_1+Z_{0}+\dots,Z_{-N+1}\} \\
&=& 
(\lim_{N\to\infty}\max\{0,\max\{0,Z_{0},\dots,Z_{0}+\dots,Z_{-N+2}\}+Z_1\})
\\
&=&
(\max\{0,
[\lim_{N\to\infty}\max\{0,Z_{0},\dots,Z_{0}+\dots,Z_{-N+2}\}]+Z_1\})
\\
&=& 
(X'_0+Z_1)^+\\
&=&X'_1.
\end{eqnarray*}
Now we prove from $n$ to $n+1$. 
\begin{eqnarray*}
X'_{n+1}&=&(X'_n+Z_{n+1})^+\\
&=&(F(T^n Z_{-\infty}^{\infty})+Z_{n+1})^+\\
&=&(\lim_{N\to\infty}
\max\{0, Z_n,Z_n+Z_{n-1},\dots,Z_n+Z_{n-1}+\dots,Z_{n-N}\}
+Z_{n+1})^+\\
&=&\lim_{N\to\infty}
\max\{0, Z_{n+1},Z_{n+1}+Z_{n},\dots,Z_{n+1}+Z_{n}+\dots,Z_{n+1-N}\}
\\
&=&
F(T^{n+1} Z_{-\infty}^{\infty}).
\end{eqnarray*}
{\bf Step 3.} Similarly to the proof of Step 1, 
\[
X_n=\max\{0,Z_n, Z_n+Z_{n-1},\dots,Z_n+\dots+Z_1,Z_n+\dots+Z_1+X_0\},
\]
and
\[
X'_n=\max\{0,Z_n, Z_n+Z_{n-1},\dots,Z_n+\dots+Z_1,Z_n+\dots+Z_1+X'_0\}.
\]
But for large $n$, both 
\[ 
Z_n+\dots + Z_1+X_0<0
\]
and 
\[
Z_n+\dots + Z_1+X'_0<0, 
\]
and so 
\[ 
X_n=X'_n=\max\{0,Z_n,
Z_n+Z_{n-1},\dots,Z_n+\dots+Z_1\}.
\]
The proof of Theorem~\ref{thm1} is complete. 

\noindent
{\bf Remark 1.} From the proof it is clear that the following extension is 
straightforward:
If $\{Z_i\}$ is coupled with a stationary and ergodic $\{Z_i'\}$, 
and $\EXP\{Z_i'\}<0$, then 
$\{X_i\}$ is coupled with a stationary and ergodic $\{X'_i\}$ such that 
\[
X'_0=\sup_{n\ge 0}V_n,
\]
where
\begin{eqnarray*} 
V_0 & = & 0,\\
V_n & = & \sum_{i=0}^{n-1}Z_{-i}',\, (n\ge 1).
\end{eqnarray*} 

\subsection{Queue length for discrete time queueing}

As an application of Theorem~\ref{thm1}
consider a discrete time queueing with constant service rate $s$, and denote
by $Y_n$ the number of arrivals in time slot $n$.
Let the initial length of the queue $Q_0$ be arbitrary non-negative integer 
valued random variable. 
Then
\[
Q_{n+1}=(Q_{n}-s+Y_{n+1})^+
\]
for $n\ge 0$. 
Put
\begin{eqnarray*} 
V_0 & = & 0,\\
V_n & = & \sum_{i=0}^{n-1}Y_{-i}-ns,\, (n\ge 1).
\end{eqnarray*} 

\begin{corollary} \label{co1}
If $\{Y_i\}$ is stationary and ergodic, and $\EXP\{Y_i\}<s$, then 
then $\{Q_i\}$ is coupled with a stationary and ergodic $\{Q'_i\}$ such that 
\[
Q'_0=\sup_{n\ge 0}V_n.
\] 
\end{corollary}

\proof
Apply Theorem~\ref{thm1} for
\[
Z_n=Y_n-s.
\]

\noindent
{\bf Remark 2.} This result together with Remark 1 has some consequences for 
network of servers (tandem of queues), when the output of a server 
$(Q_n+Y_{n+1}-Q_{n+1})$ 
is the
input of another 
server. It is easy to show that the output 
$(Q_n+Y_{n+1}-Q_{n+1})$ 
is coupled with a stationary
and ergodic sequence, and the expectations of the input and the output are
equal, so the stability condition holds for the next server, too.

\subsection{Wating time for generalized G/G/1}

This is another application of Theorem~\ref{thm1}.
According to Lindley (1952) consider the extension of the G/G/1 model.
Let $W_n$ be the waiting time of the $n$-th arrival,
$S_n$ be the service time of the $n$-th arrival, and
$T_{n+1}$ be the inter arrival time between the $(n+1)$-th and $n$-th 
arrivals.
Let $W_0$ be an arbitrary  random variable. 
Then
\[
W_{n+1}=(W_{n}-T_{n+1}+S_{n})^+
\]
for $n\ge 0$. Put
\begin{eqnarray*} 
V_0 & = & 0,\\
V_n & = & \sum_{i=0}^{n-1}(S_{-i-1}-T_{-i}),\, (n\ge 1).
\end{eqnarray*} 

\begin{corollary} 
\label{co2}
(Extension of Loynes (1964))
If $\{S_{i-1}-T_i\}$ is stationary and ergodic, 
$\EXP\{S_{i-1}\}<\EXP\{T_i\}$, then 
$\{W_i\}$ is coupled with a stationary and ergodic $\{W'_i\}$ such that 
\[
W'_0=\sup_{n\ge 0}V_n.
\] 
\end{corollary}

\proof
Apply Theorem~\ref{thm1} for
\[
Z_n=S_{n-1}-T_n.
\]

\section{Limit distribution}

In a queueing problem the properties of the limit distribution are of 
great importance. In this section we consider the special case of Section 
2.2., when the arrivals $\{Y_n\}$  are ergodic and the service rate is 
constant $s$. If $\{Y_n\}$ are weakly dependent then the tail of the limit
distribution is almost exponential, which may result in efficient algorithms 
for call admission control (cf. Duffield, Lewis, O'Connel, Russel Toomey
(1995)). The exponential tail
distribution can be derived using large deviation technique (cf. Glynn, Whitt 
(1994)).
The basic tool in this respect is the cummulant moment
generating function:
\[
\lambda(\theta)=\lim_{n \to \infty} \frac{1}{n} 
     \log \EXP \{e^{\theta \sum_{k=1}^n Y_k}\},
\]
assuming that this limit exists.
If the set
\[
\{\theta; \lambda(\theta)-\theta s < 0\},
\]
is not empty then put
\[
\delta=\sup \{\theta; \lambda(\theta)-\theta s < 0\}.
\]
Then for large $q$
\[
\PROB \{ Q>q \}\simeq e^{-\delta q}.
\]
The question is whether under the stability condition $E(Y_1)<s$ 
one has exponential tail distribution.

Glynn, Whitt  (1994) gave a positive answer under weakly dependent
$\{Y_n\}$. The limit in the definition of $\lambda(\theta)$ exists only 
under some conditions, for example, if $\{Y_n\}$ form a binary Markov
chain then $\lambda(\theta)$ can be calculated (Dembo, Zeitouni (1992)),
an explicite bound on $\PROB \{ Q>q \}$ can be given (Duffield (1994)).

For a possible extension the other problem is that $\EXP (Y_1)<s$
is only a necessary condition that the set 
$\{\theta; \lambda(\theta)-\theta s < 0\} $ is not empty, but not sufficient.
This can be seen by Jensen's inequality:
\[
\frac{1}{n} \log \EXP \{e^{\theta \sum_{k=1}^n Y_k}\}
> \frac{1}{n} \log e^{\theta \EXP \{\sum_{k=1}^n Y_k\} }
= \theta \EXP \{ Y_1\},
\]
therefore if the set $\{\theta; \lambda(\theta)-\theta s < 0\} $ is not empty
then $E(Y_1)<s$.

For long range dependent arrivals Duffield, O'Connel (1995) proved that the 
tail may not be exponential. They introduced the scaled cummulant moment
generating function:
\[
\lambda^*(\theta)=\lim_{n \to \infty} \frac{1}{v(n)} 
     \log \EXP \{e^{\theta {v(n)\over a(n)}[\sum_{k=1}^n Y_k-ns]}\},
\]
assuming that this limit exists, where $a(t)$ and $v(t)$ are 
monoton increasing functions.
If the set
\[
\{\theta; \lambda^*(\theta) < 0\},
\]
is not empty then put
\[
\delta=\sup \{\theta; \lambda^*(\theta) < 0\}.
\]
Then for large $q$
\[
\PROB \{ Q>q \}\simeq e^{-\delta v(a^{-1}(q))}.
\]
Duffield, O'Connel (1995) applied this result when $\{\sum_{k=1}^n Y_k\}$ is 
a Gaussian process with stationary increments, or Ornstein-Uhlenbeck process,
or a squared Bessel process. These examples can be motivated by multiplexing 
many sources.
For all these examples $Y_n$ is unbounded. In this section we consider bounded
$Y_n$, especially binary valued $Y_n$. Show examples such that $Q$ has larger 
tail than exponential.

\begin{proposition} \label{pro1}
There is a stationary, ergodic and binary valued $\{Y_n\}$ such that 
$\EXP \{ Y_1\}\le 1/2$ and with $s=3/4$ the queue length sequence is
stable, and for each $\delta >0$ there is $q$ such that
\[
\PROB \{ Q>q \}> e^{-\delta q}.
\]
\end{proposition}
\proof From Theorem~\ref{thm1}
\[
Q=\sup_{n\ge 0}V_n,
\]
therefore for any $n$
\[
\PROB \{ Q>q \}\ge \PROB \{ V_n>q \} = \PROB \{ \sum_{j=0}^{n-1} Y_{-j}
>ns+q \}.
\]
We show that for $\delta_i={1\over i}$, $q_i=2i^2$, $n_i=2^{i-1}$
($i>16$),
$$ 
\PROB (\sum_{j=0}^{n_i-1} Y_{-j} > {3\over 4} n_i +q_i) >2^{-\delta_i
q_i}.
$$

\noindent
{\bf Step 1.}
We present first a dynamical system given in Gy\"{o}rfi, Morvai, Yakowitz 
(1998). We will define a transformation $T$ on the unit interval.
 Consider the binary expansion 
$r_1^{\infty}$  of each real-number $r\in [0,1)$, that is, 
$r=\sum_{i=1}^{\infty} r_i 2^{-i}$.  When there are two expansions,
use the representation which contains finitely many  $1's$.
Now let 
\be
\tau(r)= \min\{i>0: r_i=1\}.
\ee
Notice that, aside from the exceptional set  $\{0\}$,  which has Lebesgue measure zero $\tau$  is finite and 
well-defined on the closed unit interval.   
The transformation is defined by 
\be
(Tr)_i=\left\{
\begin{array}{ll}
1 & \mbox{if $0<i<\tau(r)$} \\
0 & \mbox{if $i=\tau(r)$} \\
r_i & \mbox{if $i>\tau(r)$}.
\end{array}
\right. 
\ee

\noindent
{\bf Step 2.}
We show that the transformation $T$ is ergodic.\\
Notice that in fact, $Tr=r-2^{-\tau(r)}+\sum_{l=1}^{\tau(r)-1} 2^{-l}$.
  All iterations $T^k$ of $T$ for $-\infty<k<\infty$ are well defined and invertible with the exeption of the set of dyadic rationals which has Lebesgue measure zero. In the future we will neglect this set.  
Transformation $T$ could be  defined recursively as 
$$
Tr=\left\{
\begin{array}{ll}
r-0.5 & \mbox{if $0.5\le r<1$} \\
{1+T(2r) \over 2} & \mbox{if $0\le r< 0.5$.} 
\end{array}
\right.
$$
Let 
\[
S_i=\{I_0^i,\dots,I_{2^i-1}^i\}
\]
be a partition of $[0,1)$ 
where for each integer $j$ in the range $0\le j<2^i$ $I_j^i$ is defined as the set of numbers $r=\sum_{v=1}^{\infty} r_v 2^{-v}$ whose 
binary expansion $0.r_1, r_2,\dots$ starts with the bit sequence $j_1,j_2,\dots,j_i$ that is reversing the binary expansion $j_i,\dots,j_2,j_1$ of the number 
 $j=\sum_{l=1}^i 2^{l-1} j_l$. 
  Observe that in $S_i$ there are $2^i$ left-semiclosed 
intervals and each interval $I_j^i$ has length (Lebesgue measure) $2^{-i}$. 
 Now $I_{j}^i$ is mapped linearly, under $T$ onto 
$I_{j-1}^i$ for $j=1,\dots,2^i-1$. 
To confirm this, observe that  for $ j=1,\dots,2^i-1$, if $r\in I_{j}^i$ then 
\begin{eqnarray*}
Tr&=&\sum_{l=1}^{\tau(r)-1} 2^{-l}+\sum_{l=\tau(r)+1}^{\infty} r_l 2^{-l} \\
&=&r- \sum_{l=1}^i 2^{-l} ( j_l-(j-1)_l)\\
&=& \sum_{l=1}^{i} (j-1)_l 2^{-l}+\sum_{l=i+1}^{\infty} r_l 2^{-l}.  
\end{eqnarray*}
Now if $0<r\in I_0^i$ then $\tau(r)>i$ and  so $Tr\in I_{2^i-1}^i$. 
Furthermore, if $r\in I_{2^i-1}^i$ then 
$r_1=\dots=r_i=1$,  and thus conclude that $(T^{-1}r)_1=\dots=(T^{-1}r)_i=0$, that is, $T^{-1}r\in I_0^i$. Let $r\in [0,1)$ and $n\ge 1$ be arbitrary. Then   $r\in I_j^n$ for some $0\le j\le 2^n-1$. 
For all $j-(2^n-1)\le k\le j$,  
 \begin{equation} \label{iteratedshift}
T^{k}r= \sum_{l=1}^{n} (j-k)_l 2^{-l}+\sum_{l=n+1}^{\infty} r_l 2^{-l}.  
\end{equation} 
Now since $T^{-1}I_{j}^i=I_{j+1}^i$ 
for $i\ge 1$, $j=0,\dots, 2^i-2$, and the union over $i$ and $j$ of these 
sets generate the Borel $\sigma$-algebra, 
we conclude that $T$ is measurable. Similar reasoning shows that $T^{-1}$ is also measurable.  
The dynamical system $(\Omega,\F,\mu,T)$ is identified with
$\Omega=[0,1)$ and $\F$  the Borel $\sigma$-algebra on $[0,1)$, 
 $T$ being the 
transformation developed above.  Take $\mu$ to be 
 Lebesgue measure on the 
unit interval. Since transformation $T$ is measure-preserving on each set 
in the collection  $\{I_j^i : 1\le j\le 2^i-1, 1\le i<\infty\}$ 
and these intervals generate the Borel $\sigma$-algebra $\F$, 
$T$ is a stationary transformation. Now we prove that  transformation $T$ is 
ergodic as well. 
Assume $TA=A$. If $r\in A$ then $T^l r\in A$ for $-\infty<l<\infty$. 
Let $R_n: [0,1)\rightarrow \{0,1\}$ be the function $R_n(r)=r_n$. 
If $r$ is chosen uniformly on $[0,1)$ then $R_1,R_2,\dots$ is a series 
if i.i.d. random variables. 
Let $\F_n=\sigma(R_n,R_{n+1},\dots)$. By~(\ref{iteratedshift}) it is 
immediate that $A\in \cap_{n=1}^{\infty} \F_n$
and so  $A$ is a tail event. By  Kolmogorov's zero one law $\mu(A)$ is 
either zero or one.  Hence $T$ is ergodic. 
 
\noindent
{\bf Step 3.}
We define a
partition of $[0,1)$ in the following way.
Let $A_0=\emptyset$, $B_0=[0,2^{-2})$, $C_0=A_0\bigcup B_0$. 
In general, for $i\ge 1$ let 
\begin{equation} \label{defAi}
A_i=\bigcup_{j=0}^{2^{i-1}-1} T^{-j} [0,2^{-2i-1})
=\bigcup_{j=0}^{2^{i-1}-1} I_j^{2i+1}
\end{equation}
 and 
\begin{equation} \label{defBi}
B_i=\bigcup_{j=2^{i-1}}^{2^{i}-1} T^{-j} [0,2^{-2i-1})
=\bigcup_{j=2^{i-1}}^{2^{i}-1} I_j^{2i+1}.
\end{equation}
We show that
$$ \mu(A_i)=\mu(B_i)=2^{-i-2}.$$
Since $[0,2^{-2i-1})= I_0^{2i+1}$ 
it is clear for $0\le j< 2^{i}$ the sets $T^{-j} [0,2^{-2i-1})$
are disjoint. Thus 
$$ \mu(A_i)=\mu(B_i)=2^{-2i-1} 2^{i-1}.$$
 
\noindent
{\bf Step 4.} Put
\begin{equation} \label{defCi}
C_i=A_i\bigcup B_i.
\end{equation}
and 
\[
C=\bigcup_{j=0}^{\infty} C_i.
\]  
We show that 
\[
\mu(C)\le {1\over 2}.
\]
Since $A_0=\emptyset$ and for $i\ge 1$, 
\begin{eqnarray*}
A_i&=&\bigcup_{j=0}^{2^{i-1}-1} T^{-j}[0,2^{-2i-1})\\
& = & 
\left( \bigcup_{j=0}^{2^{i-1-1}-1} T^{-j} [0,2^{-2i-1}) \right) 
\bigcup
\left( \bigcup_{j=2^{i-1-1}}^{2^{i-1}-1} T^j [0,2^{-2i-1})\right)\\
&\subseteq& 
\left( \bigcup_{j=0}^{2^{i-1-1}-1} T^{-j} [0,2^{-2(i-1)-1}) \right) 
\bigcup
\left( \bigcup_{j=2^{i-1-1}}^{2^{i-1}-1} T^j [0,2^{-2(i-1)-1})\right)\\
&=& A_{i-1}\bigcup B_{i-1}\\
&=& C_{i-1}
\end{eqnarray*}
and so 
$A_i\subseteq C_{i-1}$, that is, 
\[
C=\bigcup_{i=0}^{\infty} B_i.
\]
Furthermore 
\begin{eqnarray*}
\mu(C) &\le& \sum_{i=0}^{\infty} \mu(B_i)\\
&=&\sum_{i=0}^{\infty} 2^{-2i-1} 2^{i-1}\\
&=& \sum_{i=0}^{\infty} 2^{-i-2}\\
&=&{1\over 2}.
\end{eqnarray*}  
 
\noindent
{\bf Step 5.} 
Define the binary time series $\{Y_i\}$ as
\[
Y_i(\omega)=\left\{ \begin{array}{ll}
1 & \mbox{if $T^i \omega\in C$} \\
0 & \mbox{ otherwise.}
\end{array}
\right.
\]
Clearly $\{Y_i\}$ is stationary and ergodic since the dynamical
system itself was
so. $EY_i\le 0.5$ since $\mu(C)\le 0.5$. 
If $\omega\in A_i$ then by (\ref{defAi}), (\ref{defBi}) and (\ref{defCi})  
$$\omega,T^{-1}\omega,\dots,T^{-(2^{i-1}-1)}\omega\in C$$ that is 
$$Y_0(\omega)=1,\dots,Y_{-(n_i-1)}(\omega)=1.$$
Furthermore for $i>16$,
${n_i\over 4}={2^{i-1}\over 4}> {8 i^2\over 4} = q_i$  
and so 
$n_i={3\over 4}n_i+{1\over 4}n_i>{3\over 4} n_i+q_i$. 
Thus 
$$
\PROB (\sum_{j=0}^{n_i-1} Y_{-j} > {3\over 4} n_i +q_i) \ge 
\mu(A_i).
$$
By Step 3,  for  $i>16$, 
$$
\mu(A_i)=2^{-i-2}>  2^{-2i} = 2^{-(1/i) 2 i^2}=
2^{-\delta_i q_i}.
$$ 
The proof of Proposition 1 is complete.

\bigskip

One could define 
$$
\lambda(\theta):=
\limsup_{n\to \infty} {1\over n} \log \EXP e^{\theta\sum_{i=1}^n Y_i}.
$$
The next poposition shows that for the stationary and 
ergodic time-series just defined, 
$$\{\theta ; \lambda(\theta)-\theta s < 0\}$$ 
is the emty set when $s=1$.

\begin{proposition} \label{pro2} For $\{Y_i\}$ defined in the proof of
Proposition~\ref{pro1},
$$\lim_{i\to\infty} 
{1\over n_i} \log \EXP e^{\sum_{j=0}^{n_i-1} \theta Y_j }=\theta.
$$
\end{proposition}
\proof
By Step 3 and 5 of the proof of Proposition 1
\begin{eqnarray*}
\theta&\ge& 
\limsup_{i\to\infty}{1\over n_i} \log \EXP e^{\sum_{j=0}^{n_i-1} \theta
Y_j
}\\
&\ge& 
\liminf_{i\to\infty}{1\over n_i \log_2 e} \log_2 \EXP 
2^{\sum_{j=0}^{n_i-1}
(\log_2 e)\theta Y_j }\\
&\ge& 
\liminf_{i\to\infty}{1\over n_i \log_2 e} \log_2 (2^{(\log_2 e)\theta n_i}
\mu
(A_i))\\
&=& 
\liminf_{i\to\infty} \theta+ {1\over n_i\log_2 e} \log_2 \mu (A_i)\\
&=& 
\theta+\lim_{i\to\infty}   { 2^{-i+1}\over \log_2 e} \log 2^{-i-2}=\theta.
\end{eqnarray*} 
The proof of Proposition~\ref{pro2} is complete.

\end{document}